\documentclass[a4paper]{amsart}

\RequirePackage{amsmath}
\RequirePackage{bm}
\RequirePackage{amssymb}
\RequirePackage{upref}
\RequirePackage{amsthm}
\RequirePackage{enumerate}
\RequirePackage{pb-diagram}
\RequirePackage{amsfonts}
\RequirePackage[mathscr]{eucal}
\RequirePackage{verbatim}
\RequirePackage{xr}
\RequirePackage{graphicx}
\usepackage{calc}
\usepackage{xspace}

\newcommand{\cf}{cf.\@\xspace}
\newcommand{\resp}{resp.\@\xspace}

\newcommand{\al}{\alpha}
\newcommand{\bet}{\beta}

\newcommand{\de}{\delta }
\newcommand{\e}{\epsilon}

\newcommand{\f}{\varphi}
\newcommand{\h}{\eta}

\newcommand{\ka}{\kappa}
\newcommand{\lam}{\lambda}
\newcommand{\m}{\mu}

\newcommand{\s}{\sigma}

\newcommand{\C}{\varGamma}

\newcommand{\F}{\varPhi}
\newcommand{\Lam}{\varLambda}
\newcommand{\Om}{\varOmega}

\newcommand{\Si}{\varSigma}

\newcommand{\nbdd}{\nobreakdash--}

\newcommand{\fv}[2]{#1\hspace{0pt}_{|_{#2}}}

\newcommand{\so}{{\mc S_0}}

\newcommand{\const}{\tup{const}}

\newcommand{\msp[1]}[1]{\mspace{#1mu}}

\newcommand{\R}[1][n+1]{{\protect\mathbb R}^{#1}}

\newcommand{\N}{{\protect\mathbb N}}

\newcommand{\eR}{\stackrel{\lower1ex \hbox{\rule{6.5pt}{0.5pt}}}{\msp[3]\R[]}}
\newcommand{\eN}{\stackrel{\lower1ex \hbox{\rule{6.5pt}{0.5pt}}}{\msp[1]\N}}
\newcommand{\eO}{\stackrel{\lower1ex
\hbox{\rule{6pt}{0.5pt}}}{\msc O}}

\DeclareMathOperator{\graph}{graph}

\DeclareMathOperator{\indm}{ind}

\newcommand\ra{\rightarrow}

\newcommand\hra{\hookrightarrow}

\newcommand\pa{\partial}

\newcommand{\un}{\infty}
\newcommand{\A}{\forall}

\newcommand{\set}[2]{\{\,#1\colon #2\,\}}
\newcommand{\uu}{\cup}
\newcommand{\ii}{\cap}
\newcommand{\uuu}{\bigcup}

\newcommand{\uud}{ \stackrel{\lower 1ex \hbox {.}}{\uu}}
\newcommand{\uuud}[1]{ \stackrel{\lower 1ex \hbox {.}}{\uuu_{#1}}}
\newcommand\su{\subset}

\newcommand\eS{\emptyset}
\newcommand{\sminus}[1][28]{\raise 0.#1ex\hbox{$\scriptstyle\setminus$}}

\newcommand{\wed}{\wedge}

\newcommand\ti{\times }

\newcommand{\abs}[1]{\lvert#1\rvert}

\newcommand{\norm}[1]{\lVert#1\rVert}

\newcommand{\spd}[2]{\protect\langle #1,#2\protect\rangle}

\newcommand\ch[3]{\varGamma_{#1#2}^#3}
\newcommand\cha[3]{{\bar\varGamma}_{#1#2}^#3}

\newcommand{\tit}{\textit}

\newcommand{\tup}{\textup}

\newcommand{\mc}{\protect\mathcal}
\newcommand{\msc}{\protect\mathscr}

\providecommand{\bysame}{\makebox[3em]{\hrulefill}\thinspace}

\newcommand{\cq}[1]{\glqq{#1}\grqq\,}

\newcommand{\bt}{\begin{thm}}
\newcommand{\bl}{\begin{lem}}
\newcommand{\bc}{\begin{cor}}
\newcommand{\bd}{\begin{definition}}
\newcommand{\bpp}{\begin{prop}}
\newcommand{\br}{\begin{rem}}
\newcommand{\bn}{\begin{note}}
\newcommand{\be}{\begin{ex}}
\newcommand{\bes}{\begin{exs}}
\newcommand{\bb}{\begin{example}}
\newcommand{\bbs}{\begin{examples}}
\newcommand{\ba}{\begin{axiom}}
\newcommand{\bas}{\begin{assumption}}

\newcommand{\et}{\end{thm}}
\newcommand{\el}{\end{lem}}
\newcommand{\ec}{\end{cor}}
\newcommand{\ed}{\end{definition}}
\newcommand{\epp}{\end{prop}}
\newcommand{\er}{\end{rem}}
\newcommand{\en}{\end{note}}
\newcommand{\ee}{\end{ex}}
\newcommand{\ees}{\end{exs}}
\newcommand{\eb}{\end{example}}
\newcommand{\ebs}{\end{examples}}
\newcommand{\ea}{\end{axiom}}
\newcommand{\eas}{\end{assumption}}

\newcommand{\bp}{\begin{proof}}
\newcommand{\ep}{\end{proof}}
\newcommand{\eps}{\renewcommand{\qed}{}\end{proof}}

\newcommand{\bal}{\begin{align}}

\newcommand{\bi}[1][1.]{\begin{enumerate}[\upshape #1]}
\newcommand{\bia}[1][(1)]{\begin{enumerate}[\upshape #1]}
\newcommand{\bin}[1][1]{\begin{enumerate}[\upshape\bfseries #1]}
\newcommand{\bir}[1][(i)]{\begin{enumerate}[\upshape #1]}
\newcommand{\bic}[1][(i)]{\begin{enumerate}[\upshape\hspace{2\cma}#1]}
\newcommand{\bis}[2][1.]{\begin{enumerate}[\upshape\hspace{#2\parindent}#1]}
\newcommand{\ei}{\end{enumerate}}

\newcommand\ndots{\raise 0.47ex \hbox {,}\hskip0.06em\cdots %
     \raise 0.47ex \hbox {,}\hskip0.06em} 

\newcommand{\q}{\quad}
\newcommand{\qq}{\qquad}

\newcommand{\hp}{\hphantom}

\newcommand\nd{\noindent}

\newskip\Csmallskipamount                                                
\Csmallskipamount=\smallskipamount
\newskip\Cmedskipamount
\Cmedskipamount=\medskipamount
\newskip\Cbigskipamount
\Cbigskipamount=\bigskipamount

\newcommand\cvs{\vspace\Csmallskipamount}   
\newcommand\cvm{\vspace\Cmedskipamount}

\newskip\csa
\csa=\smallskipamount

\newskip\cma
\cma=\medskipamount

\newskip\cba
\cba=\bigskipamount

\newdimen\spt
\newcommand\citem{\cvs\advance\itemno by
1{(\romannumeral\the\itemno})\hskip3pt}
\newcommand{\bitem}{\cvm\nd\advance\itemno by
1{\bf\the\itemno}\hspace{\cma}}

\newcount\itemno
\newcommand{\las}[1]{\label{S:#1}}

\newcommand{\lae}[1]{\label{E:#1}}
\newcommand{\lat}[1]{\label{T:#1}}
\newcommand{\lal}[1]{\label{L:#1}}

\newcommand{\lar}[1]{\label{R:#1}}

\newcommand{\rs}[1]{Section~\ref{S:#1}}

\newcommand{\rt}[1]{Theorem~\ref{T:#1}}
\newcommand{\rl}[1]{Lemma~\ref{L:#1}}

\newcommand{\rr}[1]{Remark~\ref{R:#1}}

\newcommand{\re}[1]{\eqref{E:#1}}
\newcommand{\frt}[1]{Theorem~\ref{T:#1} on page~\tup{\pageref{T:#1}}}
\newcommand{\frl}[1]{Lemma~\ref{L:#1} on page~\tup{\pageref{L:#1}}}

\newcommand{\fre}[1]{\eqref{E:#1} on page~\tup{\tup{\pageref{E:#1}}}}

\newskip\thmskip
\thmskip=\parindent

\newskip\hsk
\newenvironment{hinw}{\labelsep=0pt\begin{list}{}{\labelsep=0pt\itemindent=0pt\labelwidth=0pt\leftmargin=\parindent\rightmargin=0pt\partopsep=\cba}%
\item\it\nopagebreak\nopagebreak}%
{\end{list}}

\newcommand\bh{\begin{hinw}}
\newcommand{\eh}{\end{hinw}}

\newtheoremstyle{normal}
  {\cba}
  {\cba}
  {}
  {\thmskip}
  {\bfseries}
  {.}
  {\hsk}
  {}

\newtheoremstyle{abschnitt}
  {\cba}
  {\cba}
  {}
  {\thmskip}
  {\bfseries}
  {.}
  {\hsk}
  {}

\newtheoremstyle{italic}
  {\cba}
  {\cba}
  {\itshape}
  {\thmskip}
  {\bfseries}
  {.}
  {\hsk}
  {}

\newtheoremstyle{aufgaben}
  {\cba}
  {\cba}
  {}
  {}
  {\normalsize\bfseries}
  {.}
  {\hsk}
  {}

\newtheoremstyle{break}
  {\cba}
  {\cba}
  {\itshape}
  {}
  {\bfseries}
  {.}
  {\newline}
  {}

\swapnumbers
\theoremstyle{italic}
\newtheorem{thm}[subsection]{Theorem}
\newtheorem{lem}[subsection]{Lemma}
\newtheorem{prop}[subsection]{Proposition}
\newtheorem{cor}[subsection]{Corollary}

\theoremstyle{normal}
\newtheorem{rem}[subsection]{Remark}
\newtheorem{definition}[subsection]{Definition}
\newtheorem{example}[subsection]{Example}
\newtheorem{examples}[subsection]{Examples}
\newtheorem{ex}[subsection]{Exercise}
\newtheorem{note}[subsection]{}
\newtheorem{axiom}[subsection]{Axiom}
\newtheorem{assumption}[subsection]{Assumption}

\theoremstyle{aufgaben}
\newtheorem{exs}[subsection]{Exercises}

\swapnumbers

\numberwithin{equation}{section}
\numberwithin{figure}{section}

\newenvironment{textequation}[1][0.8]
{\begin{equation}
\begin{aligned}
\begin{minipage}{#1\linewidth}}
{\end{minipage}
\end{aligned}
\end{equation}
\ignorespacesafterend}

\newcommand{\btext}{\begin{textequation}}
\newcommand{\etext}{\end{textequation}}

\def\hinweis{\@startsection{subsection}{2}%
 \z@{0.7\linespacing\@plus 0.5\linespacing}{0.7\linespacing}%
{\normalfont\itshape\indent}}

\usepackage[german,english]{babel}
\usepackage{graphicx}
\RequirePackage{amsmath}
\RequirePackage{bm}
\RequirePackage{amssymb}
\RequirePackage{upref}
\RequirePackage{amsthm}
\RequirePackage{enumerate}
\RequirePackage{pb-diagram}
\RequirePackage{amsfonts}
\RequirePackage[mathscr]{eucal}
\RequirePackage{verbatim}
\RequirePackage{xr}
\RequirePackage{graphicx}
\usepackage{calc}
\usepackage{xspace}

\makeatletter
\RequirePackage{color}
\newcommand{\ann}[1]{\renewcommand{\@makefnmark}{\mbox{$^{\color{red}{\@thefnmark}}$}}%
\footnote {#1}}
\makeatother

\RequirePackage{upref}
\RequirePackage{amsthm}
\RequirePackage{enumerate}
\usepackage[mathscr]{eucal}

\usepackage{xr-hyper}

\usepackage{calc}

\newlength{\oddsidemarginlength}
\newlength{\topmarginlength}

\newcounter{numberoflines}
\newcounter{tempcc}
\oddsidemargin=\oddsidemarginlength
\evensidemargin=\oddsidemargin
\topmargin=\topmarginlength
\begin{document}

\flushbottom


\title[A general existence proof]{A general existence proof for  non-linear elliptic equations in semi-Riemannian spaces}

\author{Claus Gerhardt}
\address{Ruprecht-Karls-Universit\"at, Institut f\"ur Angewandte Mathematik,
Im Neuenheimer Feld 294, 69120 Heidelberg, Germany}
\email{gerhardt@math.uni-heidelberg.de}
\urladdr{http://www.math.uni-heidelberg.de/studinfo/gerhardt/}
\thanks{This work has been supported by the DFG}

%
\date{\today}
%


\begin{abstract}
We present a general existence proof for a wide class of non-linear elliptic equations which can be applied to problems with barrier conditions without specifying any assumptions guaranteeing the uniqueness or local uniqueness of particular solutions.

As an application we prove the existence of closed hypersurfaces with curvature prescribed in the tangent bundle of an ambient Riemannian manifold $N$ without supposing any sign condition on the sectional curvatures $K_N$. A curvature flow wouldn't work in this situation, neither the method of successive approximation.
\end{abstract}

\maketitle

\tableofcontents

\setcounter{section}{0}
\section{Introduction}

In \cite{cg:minkowski}  we considered a Minkowski type problem in $S^{n+1}$ and used for the existence proof a continuity method combined with Smale's generalization of Sard's theorem to Fredholm operators in separable Banach spaces, see \cite{smale:sard}.  

The existence proof required the usual a priori estimates and, when looking at the convex combination
\begin{equation}
t f+(1-t)f_0,
\end{equation}
where $f$ is the prescribed right-hand side and $f_0$ is one for which a \tit{unique} solution is known, which is also supposed to be a \tit{regular} point,  then  the existence of a solution for $t=1$ could be deduced by showing that either a solution for $t=1$ exists or one gets a contradiction with the uniqueness for $t=0$.

The crucial points were the uniqueness of the solution  for $t=0$ and that the operator was a local diffeomorphism near that solution which is equivalent of being a regular point in case of elliptic operators.

Now, let $N$ be a Riemannian or Lorentzian space\footnote{Since we are mainly interested in the Riemannian case, because of the particular application we have in mind, we shall use language that might not make sense in Lorentzian spaces.}, $\Om\su N$ open, connected and precompact, $F$ a symmetric, monotone and concave curvature function, and $0<f\in C^5(T(\bar\Om))$, then we consider the problem
\begin{equation}\lae{1.2}
\fv FM=f(x,\nu),
\end{equation}
where $M\su \Om$ should be a closed hypersurface of class $C^{6,\al}$ and the right-hand side is evaluated at $x\in M$ and $\nu\in T^{1,0}_x(N)$, where $\nu$ is the normal of $M$. 

We assume furthermore, that  $\Om$ is bounded by two  disjoint closed, connected and admissible hypersurfaces $M_i$, $i=1,2$, of class $C^{6,\al}$,  where $M_2$ is an upper barrier for $(F,f)$ and $M_1$ a lower barrier. Let us emphasize that the definition of \cq{upper} \resp \cq{lower} barrier also involves the direction of the continuous normal vector, i.e., the normal vector $\nu$ of $M_2$, used in the Gaussian formula, has to point into the exterior of $\Om$ and that of $M_1$ into the interior of $\Om$, see \cite[Definition 2.7.7 and Remark 2.7.8]{cg:cp} for details.

Moreover, $\bar\Om$ should be covered by a Gaussian coordinate system $(x^\al)$, $0\le \al\le n$, such that the barriers $M_i$ can be written as graphs over a closed associated hypersurface $\so$
\begin{equation}
M_i=\graph u_i=\set{x^0=u_i(x)}{x\in\so}
\end{equation}
and the $x^0$-axis is oriented  such that $u_1\le u_2$.

The solution hypersurface $M$ is also supposed to be a graph over $\so$. 

The barriers will provide a priori estimates in the $C^0$-norm. Assuming then a priori estimates in the $C^1$ and $C^2$-norms such that the curvature operator $F$ is uniformly elliptic on the solutions, we shall show that  we can find a particular solution of \re{1.2} with right-hand side $f_0$, which can be artificially forced to be unique as well as a regular point,  by defining $f_0$ appropriately, without sacrificing the a priori estimates and the barrier conditions.

Hence, the former existence proof \cite[Theorem 6.3]{cg:minkowski} with the continuity method can be applied to solve \re{1.2}. The proof works in Riemannian as well as Lorentzian spaces $N$.

As an application we generalize a previous result, \cf \cite[Theorem 3.8.1]{cg:cp}, which required that the sectional curvatures $K_N$ of the Riemannian manifold $N$ are non-positive, by dropping this restriction.

\bt\lat{1.1}
Let  $N$ be a Riemannian manifold,  $\Om\su N$ be open, connected and precompact, $F\in (K^*)$ of class $C^{5,\al}(\C_+)$, $0<f\in C^{5}(T(\bar \Om))$ and suppose that the boundary of $\Om$ has two components $M_i$, $i=1,2$, which are closed, disjoint, connected hypersurfaces of class $C^{6,\al}$,  which act as barriers for $(F,f)$ in the sense of \cite[Definition 2.7.7 and Remark 2.7.8]{cg:cp}, where $M_2$ is the upper barrier and $M_1$ the lower barrier.  Then the problem \re{1.2} has a strictly convex solution $M\su \bar\Om$ of class $C^{6,\al}$ provided $\bar\Om$ is covered by a normal Gaussian coordinate system $(x^\al)$, such that  the barriers $M_i$ can be written as graphs over some level hypersurface $\so$
\begin{equation}
M_i=\fv{\graph u_i}\so,
\end{equation}
and provided there exists a strictly convex function $\chi\in C^2(\bar\Om)$.

The solution $M$ can be written as the graph of a function $u\in C^{6,\al}(\so)$. 
\et
We emphasize that neither  curvature flows nor the method of successive approximations, that we used in \cite{cg97}, could be employed for an existence proof in this particular case.
\section{The unique particular solution}\las 2

Let $\C\su\R[n]$ be an open, convex, symmetric cone containing the positive cone and $F\in C^{5,\al}(\C)\ii C^0(\bar\C)$ be a symmetric, monotone and concave curvature function such that $\C$ is the defining cone for $F$, i.e.,
\begin{equation}
\fv F{\pa\C}=0.
\end{equation}

Notice that we do not distinguish between $F$ defined in $\C$ and $F$ defined on admissible symmetric tensors of order two for a given Riemannian metric, i.e,
\begin{equation}
F(\ka_i)=F(h_{ij})=F(h_{ij},g_{ij}),
\end{equation}
see \cite[Chapter 2.1]{cg:cp} for details.

In the Riemannian case the function $f$, which is evaluated only on unit vector fields, is certainly globally bounded. However, in the Lorentzian case, $f$ is supposed to be evaluated for unit timelike vectors and then $f$ is no longer a priori bounded. But nevertheless we can make the following assumptions on $f$ without loss of generality, see \cite[Remark 5.1.3, Remark 5.2.4 and beginning of Chapter 5.7]{cg:cp}:
\br
The function $f$ is supposed to satisfy the estimate
\begin{equation}
0<c_1\le f(x,\nu)\qq \spd{\nu}{\nu}=1,
\end{equation}
if $N$ is Riemannian, \resp
\begin{equation}
0<c_1\le f(x,\nu)\qq \spd{\nu}{\nu}=-1
\end{equation}
if $N$ is Lorentzian.

We may furthermore assume without loss of generality 
\begin{equation}\lae{2.5}
0<c_1\le f(x,\nu)\le c_2
\end{equation}
and that the barriers satisfy the corresponding inequalities strictly, i.e.,
\begin{equation}\lae{2.6}
\fv F{M_2}>f
\end{equation}
and on the set $\Si\su M_1$ of admissible points, which may be empty, there holds
\begin{equation}\lae{2.7}
\fv F{\Si}\le f-\e_1
\end{equation}
with $\e_1>0$.

These additional conditions may be assumed in the Riemannian as well as Lorentzian case without sacrificing the a priori estimates---at least for the a priori estimates that we know and used in the past.

The modifications of $f$ cited above are more general and sophisticated then we need for the present purpose, e.g., to assure \re{2.5}---with albeit different constants---we could simply replace $f$ by $\vartheta\circ f$, where $\vartheta$ is smooth and monotone satisfying
\begin{equation}
\vartheta(t)=\begin{cases}
\frac{c_1}2,&0\le t \le \frac{c_1}2,\\
t,&c_1\le t\le \frac{c_2}2,\\
c_2,&c_2\le t.
\end{cases}
\end{equation}
\er

In the Lorentzian case we assume that $N$ is globally hyperbolic with a compact Cauchy hypersurface $\so$ and that there exists a smooth global time function $x^0$. Then $N$ can be covered by a global Gaussian coordinate system $(x^\al)$, where $x^0$ is the time function and the $(x^i)$ are local coordinates for $\so$. The hypersurfaces we are interested in are all spacelike and can be written as graphs over $\so$.

The particular solution we are looking for will be a level hypersurface of $M_2$ in that part of a tubular neighbourhood  of $M_2$ which is contained in $\Om$.

First, let us establish some facts  about  tubular neighbourhoods and their foliations by level hypersurfaces.

\bl
Let $N=N^{n+1}$ be Riemannian or Lorentzian, in case $N$ is Lorentzian it is supposed to be globally hyperbolic with a compact Cauchy hypersurface $\so$, $M\su N$ a closed, oriented, spacelike\footnote{Terminology that only makes sense in a Lorentzian setting should be ignored otherwise.} hypersurface of class $C^{m,\al}$, $2\le m$, $0\le \al\le 1$, that can be written as a graph over a spacelike closed hypersurface $\so$ in a Gaussian future directed coordinate system $(x^\al)$, then there exists a tubular neighbourhood $U_{\e_0}$ of $M$ and an associated normal Gaussian coordinate system $(\tilde x^\al)$ of class $C^{m,\al}$ such that $\tilde x^0$ corresponds to the signed distance function $d=d_M$ of $M$.

The coordinate slices
\begin{equation}
M(\tau)=\{\tilde x^0=\tau\},\qq -\e_0<\tau<\e_0,
\end{equation}
are level hypersurfaces of $M=M(0)$, and if $\e_0$ is small enough, they can  also be written as graphs over $\so$
\begin{equation}\lae{2.10}
M(\tau)=\set{x^0=\f(\tau,x)}{x\in\so},
\end{equation}
where $\f$ is of class $C^{m,\al}$ in all variables such that
\begin{equation}\lae{2.11}
\dot\f>0.
\end{equation}
 Hence, there holds
 \begin{equation}\lae{2.12}
c_1\abs{\tau_1-\tau_2}\le \abs{\f(\tau_1,x)-\f(\tau_2,x)}\le c_2\abs{\tau_1-\tau_2}
\end{equation}
 with positive constants $c_1,c_2$.
 
Let $h_{ij}$ be the second fundamental form of $M(\tau)$ in the coordinate system $(x^\al)$, then $F(h_{ij})$ can also be expressed as
\begin{equation}\lae{2.13}
F(h_{ij})=F(x,\f,D\f,D^2\f),
\end{equation}
where the covariant derivatives of $\f$ are defined with respect to the metric $\s_{ij}(\f,x)$. Here, the metric of $N$ is given by
\begin{equation}\lae{2.14}
d\bar s^2=e^{2\psi}\{\s (dx^0)^2+\s_{ij}(x^0,x)dx^idx^j\}
\end{equation}
and $\s=1$, if $N$ is Riemannian, \resp $\s=-1$ in the Lorentzian case.
\el

\bp
For a proof that a tubular neighbourhood exists and that the distance function is as regular as $M$ we refer to \cite[Theorem 1.3.13]{cg:cp}.

We shall only prove \re{2.10}, \re{2.11} and \re{2.13}.

\cvm
\cq{\re{2.10}}\q Let $d$ be the signed distance function, then
\begin{equation}
M(\tau)=\{d=\tau\}.
\end{equation}

Since $M(0)=M$ is a graph over $\so$, $M=\graph u$, its normal vector $\nu$ can be expressed---apart from a sign---as 
\begin{equation}\lae{2.16}
\nu=(\nu^\al)=v^{-1}e^{-\psi}(\s,-u^i),
\end{equation}
where
\begin{equation}
v^2=1+\s\abs{Du}^2=1+\s\s^{ij}u_iu_j,
\end{equation}
\cf \cite[Chapter 1.5, Chapter 1.6]{cg:cp}.

On the other hand, $\nu$ is equal to
\begin{equation}
\nu=\pm (\bar g^{\al\bet}d_\bet),
\end{equation}
and we deduce from \re{2.14}
\begin{equation}
\frac{\pa d}{\pa x^0}=\pm v^{-1}e^\psi.
\end{equation}
A careful inspection of our choice of normal reveals that the plus sign is valid in the preceding relation
\begin{equation}
\frac{\pa d}{\pa x^0}= v^{-1}e^\psi,
\end{equation}
\cf \cite[Remark 1.5.1, Remark 1.6.1]{cg:cp}.

Hence
\begin{equation}
\frac{\pa d}{\pa x^0}>0\qq\text{in}\; M
\end{equation}
and choosing $\e_0$ small enough, this property will be valid in $U_{\e_0}$.

The implicit function theorem then yields that the hypersurfaces $M(\tau)$ can be expressed as graphs over $\so$ as in \re{2.10}, where $\f$ is of class $C^{m,\al}$ in $(\tau,x)$.

\cvm
\cq{\re{2.11}}\q Differentiating
\begin{equation}
\tau=d(\f(\tau,x),x),\qq x\in \so
\end{equation}
with respect to $\tau$ yields
\begin{equation}
1=\frac{\pa d}{\pa x^0}\dot\f,
\end{equation}
which implies \re{2.11} in view of \re{2.14}.

\cvm
\cq{\re{2.13}}\q The second fundamental form of $M(\tau)$ can be expressed as
\begin{equation}\lae{2.24}
h_{ij}v^{-1}e^{-\psi}=-\f_{ij}-\cha 000\f_i\f_j-\cha 0i0\f_j-\cha0j0\f_i-\cha ij0,
\end{equation}
where the covariant derivatives of $\f$ are those with respect to the induced metric
\begin{equation}
g_{ij}=e^{2\psi}\{\s\f_i\f_j+\s_{ij}\}.
\end{equation}

Now, $\f_{ij}$ can be expressed by $\f_{;ij}$, $(\f_{;ij})$ is the Hessian of $\f$ with respect to the metric $(\s_{ij}(\f,x))$, leading to the formula 
\begin{equation}\lae{2.26}
\begin{aligned}
e^{-\psi}v^{-1}h_{ij}=-v^{-2}\f_{;ij}+\bar h_{ij}+v^{-1}\psi_\al\tilde \nu^\al\tilde g_{ij}\end{aligned}
\end{equation}
\cf \cite[formula (2.5.11)]{cg:cp}, where $(\bar h_{ij})$ is the second fundamental form of the slices $\{x^0=\const\}$ relative to the conformal metric
\begin{equation}
\tilde g_{\al\bet}=e^{-2\psi}\bar g_{\al\bet}
\end{equation}
and where the other symbols, embellished by a tilde, are
\begin{equation}
\tilde g_{ij}=e^{-2\psi}g_{ij}
\end{equation}
and
\begin{equation}
\tilde\nu=e^\psi\nu.
\end{equation}
Hence, \re{2.13} is proved, where we apologize for the slightly ambiguous notation.
\ep

We now consider a tubular neighbourhood of $M_2$, $U_{\e_0}$, for small $\e_0$, and observe that due to our conventions
\begin{equation}
U_{\e_0}\ii\Om=\{-\e_0<d<0\}\equiv U_{\e_0}^-.
\end{equation}

We define the particular solution $M_0$ by
\begin{equation}
M_0=M(\tau_0),\qq-\e_0<\tau_0<0,
\end{equation}
where $\tau_0$ is very close to $0$. At the moment $\tau_0$ is still flexible, but it will be fixed in the uniqueness proof.

Define
\begin{equation}\lae{2.32}
f_0(x^0,x)=F(x,\f(\tau_0),\ldots)+\lam (\f(\tau_0)-x^0),
\end{equation}
where $0<\lam$ is very large.

Notice that $\tau_0=\tau_0(\lam)$ should always be chosen such that
\begin{equation}
\begin{aligned}
\tfrac12 F(x,\f(\tau_0),\ldots)&\le F(x,\f(\tau_0),\ldots)+\lam(\f(\tau_0)-\f(0))\\
&\le F(x,\f(0),\ldots).
\end{aligned}
\end{equation}
The last inequality is automatically satisfied, if $\lam$ is large, since $\f(\tau_0)<\f(0)$.

These assumptions imply that the barriers will also be barriers for the combinations
\begin{equation}\lae{2.34}
tf+(1-t)f_0,\qq-\de\le t\le 1+\de,
\end{equation}
for small $0<\de$.

We can now prove the uniqueness of $M(\tau_0)$ and consider first the Riemannian case.

\bl\lal{2.3}
Let $N$ be Riemannian. Then $M(\tau_0)$ is the unique solution of
\begin{equation}
\fv F{M(\tau_0)}=f_0
\end{equation}
among all admissible hypersurfaces $M\su\Om$ which can be written as graphs over $\so$ in the coordinate system $(x^\al)$, if $\lam$ is large, $\lam\ge \lam_0$, where $\lam_0=\lam_0(M_2,\Om)$.
\el
\bp
Let $M=\graph u$ be another solution. The geometric quantities of $M$ will be denoted by $h_{ij}, g_{ij}$, etc..

We distinguish two cases.
\bh
Case $1$
\eh
Suppose that
\begin{equation}
\sup_\so(\f(\tau_0)-u)>0
\end{equation}
and let $x_0\in\so$ be a point where the supremum is realized. Then
\begin{equation}
\f_{,ij}\le u_{,ij}\q\wed\q \f_i=u_i
\end{equation}
where a comma indicates partial derivatives and where we simply write $\f$ instead of $\f(\tau_0)$.

From \re{2.26} we then deduce
\begin{equation}\lae{2.38}
\begin{aligned}
h_{ij}v^{-1}e^{-\psi}&\le -v^{-2}\f_{,ij} +v^{-2}\ch ijk(u)\f_k+\bar h_{ij}\\
&\hp{=}\;+v^{-1}\psi_\al\tilde\nu^\al\tilde g_{ij}
\end{aligned}
\end{equation}
which implies
\begin{equation}\lae{2.39}
h_{ij}\le \tilde h_{ij}+c_{ij}\equiv b_{ij},
\end{equation}
where $\tilde h_{ij}$ is the second fundamental form of $\f\equiv\f(\tau_0)$ and where the tensor $c_{ij}$ depends on $(u-\f)$ such that
\begin{equation}
\norm{c_{ij}}\le c \,\abs{u-\f}
\end{equation}
with a uniform constant $c$. Notice that this is a pointwise estimate, i.e., presently in $x_0$. The metric on the left-hand side could be $g_{ij}(\f), \s_{ij}(\f)$ or $\s_{ij}(u)$.

Let $\tilde \ka_i$ be the eigenvalues of $(b_{ij})$ with respect to $g_{ij}$, then
\begin{equation}
\ka_i\le\tilde\ka_i,
\end{equation}
where all eigenvalues are labelled such that
\begin{equation}
\ka_1\le\ka_2\le\cdots\le\ka_n,
\end{equation}
etc., and where we observe that  the pair $(b_{ij}, g_{ij})$ is admissible in $x_0$, \cf \cite[Lemma 2.7.3]{cg:cp}.

Hence we infer, in view of the monotonicity and concavity of $F$,
\begin{equation}\lae{2.43}
F(\ka_i)\le F(\tilde\ka_i)\le F(\bar\ka_i)+\sum_iF_i(\bar\ka_i)(\tilde\ka_i-\bar\ka_i),
\end{equation}
where $\bar\ka_i$ are the principal curvatures of $M(\tau_0)$ in $(\f(\tau_0),x_0)$. 

Here, we also used the fact that $\C$ is convex.

Let $\hat\ka_i$ be the eigenvalues of $\tilde h_{ij}$ with respect to the metric $g_{ij}$, then we deduce from
\begin{equation}
b_{ij}=\tilde h_{ij}+c_{ij}\le \tilde h_{ij}+c\,\abs{u-\f}g_{ij}
\end{equation}
that
\begin{equation}
\tilde \ka_i\le \hat\ka_i+c\, \abs{u-\f},
\end{equation}
\cf \cite[Lemma 2.7.3]{cg:cp}. 

Hence, we only have to estimate
\begin{equation}
\hat\ka_i-\bar\ka_i
\end{equation}
from above.

To compare $\hat\ka_i$ and $\bar\ka_i$ we use the Courant-Fischer-Weyl maximum-mini\-mum principle, which says that the $i$-th eigenvalue $\ka_i$ of a symmetric matrix $A$, in the above ordering,  is determined by
\begin{equation}
\ka_i=\max \set{d(E)}{\dim E\le i-1},
\end{equation}
where $E\su \R[n]$ is a subspace and
\begin{equation}
d(E)=\min\set{\spd {A\xi}\xi}{\xi\in E^\perp,\; \abs \xi=1},
\end{equation}
\cf \cite[p. 26--29]{courant-hilbert-I}.

In the present situation we have the same covariant tensor $\tilde h_{ij}$ but different metrics
\begin{equation}
g_{ij}=e^{2\psi(u)}\{\f_i\f_j+\s_{ij}(u)\}
\end{equation}
and
\begin{equation}
\tilde g_{ij}\equiv g_{ij}(\f)=e^{2\psi(\f)}\{\f_i\f_j+\s_{ij}(\f)\}.
\end{equation}

Let $E$ denote a subspace of $T^{0,1}_{p_0}(M_0)$, $p_0=(\f(\tau_0,x_0),x_0)$, and let  $\xi\in T^{1,0}_{p_0}(M_0)$, then we can define
\begin{equation}
E^\perp=\set{\xi\in T^{1,0}_{p_0}(M_0)}{\h_i\xi^i=0\q\A\,\h\in E}
\end{equation}
independent of any metric, and the maximum-minimum principle for the pair $(\tilde h_{ij}, g_{ij})$ can be rephrased as
\begin{equation}
d(E)=\min\set{\tilde h_{ij}\xi^i\xi^j}{\xi\in E^\perp, \; g_{ij}\xi^i\xi^j=1}.
\end{equation}

Let $\tilde d(E)$ be the corresponding value for the pair $(\tilde h_{ij},\tilde g_{ij})$ and consider an arbitrary $0\ne \xi$, then
\begin{equation}\lae{2.53}
\begin{aligned}
\frac{\tilde h_{ij}\xi^i\xi^j}{g(\xi,\xi)}-\frac{\tilde h_{ij}\xi^i\xi^j}{\tilde g(\xi,\xi)}= \frac{\tilde h_{ij}\xi^i\xi^j\{\tilde g(\xi,\xi)-g(\xi,\xi)\}}{g(\xi,\xi)\tilde g(\xi,\xi)}\le c\,\abs{u-\f},
\end{aligned}
\end{equation}
where $c=c(M_2,\Om)$ as one easily checks.

Hence, we deduce
\begin{equation}
d(E)\le \tilde d(E)+c\,\abs{u-\f}
\end{equation}
yielding
\begin{equation}
\hat\ka_i\le \bar\ka_i+c\,\abs{u-\f}.
\end{equation}

Inserting these estimates in \re{2.43} we obtain
\begin{equation}
F(x_0,\f,\ldots)+\lam(\f-u)\le F(x_0,\f,\dots)+c(\f-u)
\end{equation}
where $c=c(M_2,\Om)$; a contradiction, if $\lam>c$.

\bh
Case $2$
\eh

Suppose that
\begin{equation}\lae{2.57}
\sup_\so(\f-u)\le 0,
\end{equation}
then
\begin{equation}
M\su U_{\e_o}^-.
\end{equation}

Let $p_0\in M$ be a point such that
\begin{equation}
\tau=d(p_0)=\sup_Md,
\end{equation}
where $d=d_M$, then $p_0=(\f(\tau,x_0),x_0)$ and
\begin{equation}\lae{2.60}
h_{ij}\ge \bar h_{ij}
\end{equation}
in $p_0$, where $h_{ij}$ is the second fundamental form of $M$ and $\bar h_{ij}$ the second fundamental form of $M(\tau)$.

Inequality \re{2.60} follows from \re{2.24}, since $M(\tau)$ touches $M$ from above, i.e.,
\begin{equation}
u_{ij}\le \f_{ij}
\end{equation}
and all other terms agree on the right-hand side of \re{2.24} when $u$ is replaced by $\f$.

Thus we deduce
\begin{equation}
F(x_0,\f(\tau_0),\ldots)+\lam (\f(\tau_0)-\f(\tau))\ge F(x_0,\f(\tau),\ldots)
\end{equation}
from which we further conclude with the help of \re{2.12}
\begin{equation}
\begin{aligned}
-c_1\lam\,\abs{\tau-\tau_0}&\ge\lam (\f(\tau_0)-\f(\tau))\\
&\ge F(x_0,\f(\tau),\ldots)-F(x_0,\f(\tau_0),\ldots)\\
&\ge -c\,\abs{\tau-\tau_0};
\end{aligned}
\end{equation}
a contradiction, if
\begin{equation}
\lam>\frac c{c_1},
\end{equation}
since obviously $\tau\ne\tau_0$, for otherwise $M=M(\tau_0)$ because of \re{2.57}.
\ep

Let us now consider the Lorentzian case. The first part of the proof of the previous lemma has then to be modified. However, even the modification will only work, if
\begin{equation}\lae{2.65}
v^{-1}(x_0)
\end{equation}
can be uniformly bounded independently of $\lam$.

Let us recall that
\begin{equation}\lae{2.66}
v^2=1-\s^{ij}u_iu_j=1-\s^{ij}(u)\f_i\f_j
\end{equation}
in the point $x_0$.

\br\lar{2.4}
If $M_2$ would be a coordinate slice $\{x^0=\const\}$, then $D\f(0,x)$ would vanish, and hence, by choosing $\tau_0$ small, we could also guarantee that $v^{-1}(x_0)$ would be bounded.

In \cite[Proposition 3.2]{br:mean} it is proved that there exists a time function $\tilde x^0$ such that
\begin{equation}
M_2=\{\tilde x^0=0\}.
\end{equation}
In that paper $M_2$ was assumed to be smooth, however, if $M_2$ is of class $C^{m,\al}$, then the new time function is also of class $C^{m,\al}$ and from \cite[Theorem 1.4.2]{cg:cp} we then deduce that there exists a corresponding Gaussian coordinate system $(\tilde x^\al)$ of class $C^{m-1,\al}$. 

Thus, assuming $M_2$ to be of class $C^{8,\al}$, then the new  coordinate system is of class $C^{7,\al}$, $M_2$ is a coordinate slice of class $C^{7,\al}$ and the new Christoffel  symbols are class $C^{5,\al}$, which is necessary in the next section.

Moreover, any closed, connected, spacelike hypersurface, that could be written as a graph over $\so$ in the old coordinate system, is also a graph in the new coordinate system over $M_2$, \cf \cite[Proposition 1.6.3]{cg:cp}.
\er
\bl
Let $N$ be Lorentzian with compact Cauchy hypersurface $\so$, then $M(\tau_0)$ is the unique solution of 
\begin{equation}
\fv F{M(\tau_0)}=f_0
\end{equation}
among all admissible spacelike hypersurfaces $M\su\Om$,  if $\lam$ is large, $\lam\ge \lam_0$, where $\lam_0=\lam_0(M_2,\Om)$ and $\abs{\tau_0}$ small.
\el
\bp
In view of inequality \re{2.6} we may, without loss of generality, assume  that the upper barrier $M_2$ is as smooth as $\so$ allows, i.e., we may assume $M_2\in C^{8,\al}(\so)$ at least. Then, following \rr{2.4}, we switch to a new Gaussian coordinate system $(\tilde x^\al)$ such that $M_2$ is a coordinate slice; in order to avoid any confusion with the notations in the Riemannian case,    we continue to denote the coordinates by $(x^\al)$; please, observe that this coordinate system will be the standard coordinate system in the Lorentzian case for the rest of the paper.

The only point where we need a new proof is \tit{Case} $1$, i.e.,
\begin{equation}
\sup_\so(\f(\tau_0)-u)>0.
\end{equation}

The previous inequalities \re{2.38} and \re{2.39} remain valid with the exception that the tensor $c_{ij}$ can now be estimated by 
\begin{equation}
\norm{c_{ij}}\le c(\abs{v-\bar v}+\abs{u-\f})
\end{equation}
where $v$ and $\bar v$ are the corresponding quantities for $u$ \resp $\f$.

From
\begin{equation}
(v-\bar v)(v+\bar v)=v^2-\bar v^2=\{\s^{ij}(\f)-\s^{ij}(u)\}\f_i\f_j
\end{equation}
we deduce
\begin{equation}
\abs{v-\bar v}\le c\, \abs{u-\f},
\end{equation}
where $c$ is independent of $\lam$.

The remaining inequality that has to be checked is  inequality \re{2.53}. The crucial term is $g(\xi,\xi)$, which can be estimated from below by
\begin{equation}
g(\xi,\xi)=g_{ij}\xi^i\xi^j\ge v^2 \s_{ij}(u)\xi^i\xi^j,
\end{equation}
since $v^2$ is the smallest eigenvalue of $g_{ij}$ with respect to the metric $\s_{ij}(u)$.
\ep

\br
The proofs of the uniqueness lemmata show that only the upper barrier $M_2$ is necessary to define a unique particular solution as long as the potential solution hypersurfaces are all contained in $\Om$.

Moreover, one cannot use the lower barrier $M_1$ in a similar fashion to define a particular solution as one easily checks.
\er

Let us finally prove that $\f(\tau_0)$ is a regular point for the associated Fredholm operator.

Let $\h\in H^{5,p}(\so)$, $n<p<\un$, and $\e$ small. Then $M_\e=\graph u_\e$ with
\begin{equation}
u_\e=\f(\tau_0)+\e\h
\end{equation}
is an admissible hypersurface with second fundamental form $h_{ij}(\e)$. 

Consider the operator
\begin{equation}
G(u_\e)=F(x,u_\e,Du_\e,D^2u_\e)-f_0(u_\e,x),
\end{equation}
\cf \re{2.13}.  $G$ is a Fredholm operator defined in a small neighbourhood of $\f(\tau_0)$ in $H^{5,p}(\so)$ with target space $H^{3,p}(\so)$ such that
\begin{equation}
\indm G=0.
\end{equation}
hence, $\f(\tau_0)$ is a regular point for $G$ if and only if
\begin{equation}
N(DG(\f(\tau_0)))=\{0\}.
\end{equation}

But this is satisfied for large $\lam$, since
\begin{equation}
\begin{aligned}
\spd{DG}{\h}=\frac{d}{d\e}\fv{G(u_\e)}{\e=0}=-a^{ij}\h_{;ij}+b^i\h_i+c\h+\lam\h,
\end{aligned}
\end{equation}
where $c$ is uniformly bounded and $a^{ij}$ uniformly positive definite. Hence, $N(DG(\f(\tau_0)))=\{0\}$, if $\lam>c$.

\section{The existence theorem}\las 3

The existence is proved by a continuity method using Smale's infinite dimensional version of Sard's theorem  \cite{smale:sard}. We used this method in \cite[Section 6]{cg:minkowski} to solve the Minkowski problem in the sphere and the proof given there can be carried over almost directly. However, for the convenience of the reader we shall repeat the proof in the present setting.

\br
To simplify the presentation we shall only treat the Riemannian case; the proof will also work in a Lorentzian setting, only the terminology would have to be slightly adapted.
\er

Consider the Banach spaces $E_1$, $E_2$ defined by
\begin{equation}
E_1=H^{5,p}(\so)
\end{equation}
and
\begin{equation}
E_2=H^{3,p}(\so)
\end{equation}
for some fixed $n<p<\un$, such that $H^{m,p}(\so)\hra C^{m-1,\al}(\so)$. 

Let $\tilde\Om\su E_1$ be an open set such that $u\in\tilde\Om$ implies  $M(u)=\graph u$ is an admissible hypersurface for the open, symmetric, convex cone $\C\su\R[n]$ and contained in $\Om$. We then define
\begin{equation}
\F:\tilde\Om\ra E_2
\end{equation}
by
\begin{equation}
\F(u)=F(h_{ij})-f(x,u,Du)=F(x,u,Du,D^2u)-f(x,u,Du),
\end{equation}
\cf  \re{2.13} and \re{2.16}. 

All possible solutions of $\F=0$ are strictly contained in $\tilde\Om$, if $\tilde\Om$ is specified by the requirements
\begin{equation}\lae{3.5}
u_1<u<u_2,
\end{equation}
where $u_i$ are the barriers,
\begin{equation}\lae{3.6}
\abs{Du}^2=\s^{ij}u_iu_j<c,
\end{equation}
where $c$ is a large constant, and
\begin{equation}\lae{3.7}
\ka_i<\bar\ka,
\end{equation}
where $\ka_i$ are the principal curvatures of $\graph u$,  and $\bar\ka$ large.

The constants should all be chosen such that  a solution of $\F=0$ satisfies these estimates strictly. Notice that \re{3.5} follows immediately from the maximum principle because of \re{2.5} and \fre{2.6}. 

It is well known that:

\bl
$\F$ is a nonlinear Fredholm operator of index zero.
\el

Recall that $w\in E_2$ is said to be a \tit{regular} value for $\F$, if either $w\notin  R(\F)$, or if for any $u\in \F^{-1}(w)$ $D\F(u)$ is surjective.

Smale \cite{smale:sard} proved that for separable Banach spaces $E_i$ and for Fredholm maps $\F$ the set of regular values in $E_2$ is  dense, if $\F$ is of class $C^k$ such that
\begin{equation}
k>\max(\tup{ind\,}\F,0).
\end{equation}
All requirements are satisfied in the present situation.

Next we consider the combination in \fre{2.34}
and assume furthermore that the constants used in the definition of $\tilde\Om$ are such that all possible solutions of
\begin{equation}
F=tf+(1-t)f_0,\qq -\de\le t\le 1+\de,
\end{equation}
in $\Om$ satisfy the corresponding estimates strictly. We also call the particular solution, the existence of which is proved in \rs{2}, $u_0$. The symbol $\f$ will have a different meaning in the following.

Define
\begin{equation}
\Lam: \tilde\Om\ti (-\de,1+\de)\ra E_2
\end{equation}
 by
 \begin{equation}
\Lam (u,t)=F(h_{ij})-(tf+(1-t)f_0).
\end{equation}
Then $\Lam$ is also a Fredholm operator such that $\tup{ind\,}\Lam(\cdot,t)=0$ for fixed $t$, and, if $w\in E_2$ is a regular value for $\Lam$, then
\begin{equation}\lae{6.37}
\tup{ind\,}\Lam=1\qq\A\, (u,t)\in\Lam^{-1}(w).
\end{equation}

Recall that
\begin{equation}
\tup{ind\,}\Lam =\dim N(D\Lam)-\dim \tup{coker\,} (D\Lam).
\end{equation}

The relation \re{6.37} will be proved in \rl{6.5} below.
\bt\lat{3.3} 
Let $N=N^{n+1}$ be Riemannian, $\Om\su N$ open connected and precompact and assume that $F$ is a symmetric, monotone and concave curvature function such that $F\in C^5(\C)\ii C^0(\bar\C)$ and let $0<f\in C^5(T^{1,0}(\bar\Om))$. $\Om$ should be covered by a Gaussian coordinate system $(x^\al)$ with associated closed, connected hypersurface $\so$.
We furthermore assume that $\pa\Om$ has two boundary components $M_i$ which are closed hypersurfaces of class $C^{6,\al}$ which act as barriers for $(F,f)$ and can be written as graphs in the coordinate system $(x^\al)$.  Moreover, for all possible solutions of 
\begin{equation}
\fv FM=\tilde f,
\end{equation}
with $M\su\Om$, that are graphs over $\so$,  where $\tilde f$ has the same properties as $f$ and satisfies the same structural conditions as $f$, especially \re{2.5}, \re{2.6} and \fre{2.7}, uniform a priori estimates of the form \re{3.5}, \re{3.6} and \re{3.7} are valid. Then the problem
\begin{equation}
\fv FM= f,
\end{equation}
has an admissible solution $M\su\Om$ of class $C^{6,\al}$.
\et

\bp
Consider the Fredholm map $\Lam =\Lam(u,t)$. The theorem will be proved, if we can show that there exists $u\in\Om$ such that
\begin{equation}
\Lam(u,1)=0.
\end{equation}
Note, that once we have a solution of class $H^{5,p}$ the embedding theorem and the Schauder estimates will provide the final regularity of the solution.

On the other hand, as we have proved in \frl{2.3}, there exists a unique solution of the equation
\begin{equation}\lae{6.47}
\Lam(u,0)=0,
\end{equation}
namely, $u=u_0$, the particular solution, which is also a regular point for $\Lam(\cdot,0)$, or equivalently, $(u_0,0)$ is a regular point for $\Lam$. 

Without loss of generality we may assume $0\notin R(\Lam(\cdot,1))$, for otherwise we have nothing to prove, and thus, $0$ is also regular value for $\Lam(\cdot,1)$.

Let $\e>0$ be small, then there exists a
\begin{equation}
w_\e\in B_\e(0)\su E_2,
\end{equation}
such that
\begin{equation}
tf+(1-t)f_0+w_\e>0\qq\A\,-\de\le t\le 1+\de,
\end{equation}
$w_\e\in R(\Lam(\cdot,0))$, and such that  $w_\e$ is a regular value for $\Lam(\cdot,0)$, $\Lam(\cdot,1)$ and $\Lam$.
The fact that $w_\e$ is also a regular value for $\Lam(\cdot,1)$ is due to 
\begin{equation}
\Lam(\cdot,1)^{-1}(w_\e)=\eS,
\end{equation}
\cf the reasoning below.

Set
\begin{equation}
\C_\e=\Lam^{-1}(w_\e),
\end{equation}
then $\C_\e\ne\eS$ and $\C_\e$ is a $1$-dimensional submanifold without boundary.

The intersection
\begin{equation}
\tilde\C_\e=\C_\e\ii (E_1\ti [0,1])
\end{equation}
is then compact, because of the a priori estimates, and it consists of finitely many closed curves or segments.

We want to prove that there is $u_\e\in \tilde\Om$ such that $(u_\e,1)\in \tilde \C_\e$. Suppose this were not the case, then consider a point $(\bar u_\e,0)\in\tilde\C_\e$. Such points exist by assumption. Moreover, the $1$-dimensional connected submanifold $M_\e\su \C_\e$ containing $(\bar u_\e,0)$ can be expressed near $(\bar u_\e,0)$ by
\begin{equation}\lae{6.52}
M_\e=\set{(\f(t),t)}{-\de<t<\de},
\end{equation}
where $\f\in C^1$, $\f(0)=\bar u_\e$, and
\begin{equation}
\Lam(\f(t),t)=w_\e,
\end{equation}
since by assumption $D_1\Lam(\bar u_\e,0)$ is an isomorphism and the implicit function theorem can be applied.

Let $\tilde M_\e\su M_\e\ii\tilde \C_\e$ be a connected component containing $(\bar u_\e,0)$, then $\tilde M_\e$ isn't closed because of \re{6.52}, and hence has two endpoints, see \cite[Appendix]{milnor}. One of them is $(\bar u_\e,0)$ and the other also belongs to $\Lam(\cdot,0)^{-1}(w_\e)$ and can therefore be expressed as
\begin{equation}
(\tilde u_\e,0),
\end{equation}
where $\tilde u_\e\ne \bar u_\e$ because of the implicit function theorem.

Hence we have proved that the assumption
\begin{equation}
\Lam(\cdot,1)^{-1}(w_\e)=\eS
\end{equation}
implies
\begin{equation}
\#\Lam(\cdot,0)^{-1}(w_\e) >1.
\end{equation}

However, we shall show that $\Lam(\cdot,0)^{-1}(w_\e)$ contains only one point, if $\e$ is small.

Indeed, let $\bar u_\e\in\Lam(\cdot,0)^{-1}(w_\e)$, then the $\bar u_\e$ converge to the unique solution $u_0$ of \re{6.47}. Thus, if $\e$ is small, all $\bar u_\e$ are contained in an open ball
\begin{equation}
B_\rho(u_0)\su \tilde\Om,
\end{equation}
where $\F=\Lam(\cdot,0)$ is a diffeomorphism due to the results at the end of \rs{2}, hence there exists just one solution of the equation 
\begin{equation}
\Lam(\bar u_\e,0)=w_\e.
\end{equation}

Thus we have proved that there exists a sequence
\begin{equation}
u_\e\in \Lam(\cdot,1)^{-1}(w_\e),
\end{equation}
if $\e$ tends to zero. A subsequence will then converge to a solution $u$ of
\begin{equation}
\Lam(u,1)=0.\qedhere
\end{equation}
\ep

It remains to prove the following lemma:
\bl\lal{6.5}
Let $\Lam$ be defined as above, then
\begin{equation}
\tup{ind\,}\Lam=1.
\end{equation}
\el

\bp
Let $(u_0,t_0)\in \tilde\Om\times (-\de,1+\de)$ be an \tit{arbitrary} point, where we may assume that $t_0=1$, since $\indm\Lam$ is continuous. 

We distinguish two cases:

\bh
\tit{Case} $1$\tup{:} \q\; $(f-f_0)\in R(D\F(u_0))$
\eh

We have
\begin{equation}
D\Lam=(D_1\Lam,-(f-f_0)),
\end{equation}
where all derivatives are evaluated at $(u_0,1)$ \resp $u_0$.  Then we deduce
\begin{equation}
\dim N(D\Lam)=\dim N(D_1\Lam)+1=\dim N(D\F)+1,
\end{equation}
for let 
\begin{equation}
D_1\Lam u_1=f-f_0,
\end{equation}
then
\begin{equation}
N(D\Lam)=N(D\F)\times \{0\}\oplus \langle{(u_1,1)}\rangle
\end{equation}
as one easily checks, and of course there holds
\begin{equation}
R(D\Lam)=R(D\F).
\end{equation}

Notice that this argument is also valid, if 
\begin{equation}
f(u_0,\cdot)=f_0(u_0,\cdot).
\end{equation}

\bh
\tit{Case} $2$\tup{:}\q\; $(f-f_0)\notin R(D\F(u_0))$
\eh

In this case
\begin{equation}
R(D\Lam)=R(D_1\Lam)\oplus \langle{(f-f_0)}\rangle
\end{equation}
and
\begin{equation}
N(D\Lam)=N(D_1\Lam)\times \{0\},
\end{equation}
hence
\begin{equation}
\tup{ind\,}\Lam=\tup{ind\,} \F+1=1
\end{equation}
in both cases.
\ep

\section{Proof of \rt{1.1}}\las 4
The theorem has already been proved in Riemannian manifolds $N$ the sectional curvatures of which satisfy $K_N\le 0$, \cf  \cite[Theorem 3.8.1]{cg:cp}, with the help of a curvature flow.

Now we want to apply \frt{3.3}. The barriers ensure that the hypersurfaces stay in $\Om$. Since we consider convex hypersurfaces in a normal Gaussian coordinate system uniform $C^1$-estimates are valid. Hence, it remains to prove that the principal curvatures of all solutions of
\begin{equation}
\fv FM=f
\end{equation}
are uniformly bounded from above. Since $f\ge c>0$ and $F$ vanishes on $\pa\C_+$, the principal curvatures then stay in a uniformly compact subset of $\C_+$ and $F$ will be uniformly elliptic.

To prove the upper bound for $\ka_i$ we argue as in the proof of \cite[Theorem 3.8.1]{cg:cp}, however, with a minor modification, since we have to employ a strictly convex function $\chi$ explicitly, which could be hidden in case $K_N\le0$.

For better compatibility with the former result, look at the equivalent equation
\begin{equation}
\fv{\F(F)}M=\F(f)\equiv \tilde f
\end{equation}
where $\F(t)=\log t$.

The second fundamental form satisfies an elliptic equation to which we want apply the maximum principle. 

Let the functions 
$\f$ and $w$ be  defined respectively  by
\begin{align}
\f&=\sup\set{{h_{ij}\h^i\h^j}}{{\norm\h=1}},\\
w&=\log\f+\lambda \log v+\m\chi,\lae{3.1.6}
\end{align}
where $\lambda,\m$ are large positive parameters. 

As we proved in \cite[Lemma 3.8.3]{cg:cp}
$w$ is a priori bounded for a suitable choice of $\lambda,\m$. The only difference to the present definition of $w$ is that instead of $\chi$ we used the function $u$ defined on $M$. However, $\chi$ in the general case  behaves similar to $u$ in the special case $K_N\le0$.

Notice also that the former equation for the second fundamental form was a parabolic equation, but its elliptic version can immediately be recovered.

Thus, all prerequisites are in place to apply \frt{3.3}.  


\bibliographystyle{hamsplain}
\providecommand{\bysame}{\leavevmode\hbox to3em{\hrulefill}\thinspace}
\providecommand{\href}[2]{#2}



\end{document}